\documentclass[11pt, reqno]{amsart}
\usepackage{amsmath, amsthm, amscd, amsfonts, amssymb, graphicx, color}
\usepackage[bookmarksnumbered, colorlinks, plainpages]{hyperref}

\makeatletter \oddsidemargin.9375in \evensidemargin \oddsidemargin
\DeclareMathOperator{\Rad}{Rad}
\DeclareMathOperator{\Lip}{Lip}
\newcommand{\tnorm}[1]{{\left\vert\kern-0.25ex\left\vert\kern-0.25ex\left\vert #1 
		\right\vert\kern-0.25ex\right\vert\kern-0.25ex\right\vert}}

\newtheorem{thm}{Theorem}

\newtheorem{cor}[thm]{Corollary}
\newtheorem{lem}[thm]{Lemma}

\theoremstyle{thmstyletwo}%
\newtheorem{ex}{Example}%
\newtheorem{rem}{Remark}%

\theoremstyle{thmstylethree}%

\begin{document}
\setcounter{page}{1}
	
\title{\bf BSE-Properties of Second Dual of Banach Algebras}

\author{Maryam Aghakoochaki and Ali {Rejali}$^*$}

\thanks{2020 Mathematics Subject Classifcation. D8, H51}
\thanks{*Corresponding auther}

	
\begin{abstract}
 Let $A$ be a commutative semisimple Arens regular unital Banach algebra.
The correlation between the BSE-property of the Banach algebra $A$ and its second duals are assessed. It is found and approved that if $A$ is a BSE-algebra, then so is $A^{**}$.
The opposite correlation will hold in certain conditions.  The correlation of the
 BSE-norm property of the Banach algebra and its second dual are assessed and examined. 
It is revealed that,
 if $A$ is a commutative Arens regular unital Banach algebra
where   $A^{**}$  is semisimple, then, $A$ and $A^{**}$ are BSE-norm algebra.
\end{abstract}

\maketitle

\section{Introduction}\label{sec1}

          Let $A$ be a Banach algebra. The researcher in \cite{Are} and  \cite{Are2}, introduced the two Arens multiplications on $A^{**}$ which convert it into a Banach algebra. A Banach algebra is an Arens regular if the Arens multiplications coincide with its second dual.

The  Bochner-Schoenberg-Eberlein (BSE) is derived from the famous theorem proved in 1980 by Bochner and Schoenberg for the group of real numbers; \cite{SB} and \cite{P1}. The researcher in \cite{N2},  revealed that if $G$ is any locally compact abelian group, then the group algebra $L_{1}(G)$ is a BSE algebra.
 The researcher in \cite{P1},\cite{E7},\cite{E8} assessed the commutative Banach algebras that meet the Bochner-Schoenberg-Eberlein- type theorem and explained their properties.
They are introduced and assessed in    \cite{E6} the first and second types of BSE algebras. This concept is  expanded in \cite{F1} and \cite{K1}.
 The concept of an n-BSE function and Banach algebras $C_{\textup{BSE}(n)}(\Delta(A))$, for each $n\in\mathbb N$ is introduced in  \cite{E7},  where the
  correlation between these concepts and the BSE property is assessed.
In 1992, a class of commutative Banach algebras was introduced by the same researchers in the form of the following equation:   
  $$
  \widehat{A}= C_{\textup{BSE}}(\Delta(A))
  $$
    They proved that the commutative uniform  $C^{*}$-algebras belong to this category; \cite{E7}.

 The concept of BSE-norm algebra, where $L_{1}(G)$ is a BSE-norm algebra, for any locally compact group, Hausdorff, and abelian group $G$ is devised in \cite{E8}.
Another essential reference for BSE-norm algebras is \cite{HGDU}.
The BSE property for some other Banach algebras, including the direct sum of Banach algebras is assessed in \cite{KL}. 
 That the Lipschitz algebra $\Lip_{\alpha}(K, A)$ is a BSE-algebra if and only if $A$ is a BSE-algebra, where $K$ is a compact metric space, $A$ is a commutative unital semisimple Banach algebra, and $0<\alpha\leq 1$ is proved in  \cite{AK}.

In this article, $A$ is a commutative unital Banach algebra, where $A$ is Arens regular and  $A^{**}$ is semisimple. 
That the  $C_{\textup{BSE}}(\Delta(A^{**}))\cong C_{\textup{BSE}}(\Delta(A))$ will be proved in the article first, followed by proving that if $A$ is a BSE-algebra, so is $A^{**}$.
 It  will be concluded  that if $A^{**}$ is a BSE-algebra and  $\widehat{A}\mid_{\Delta(A)}\cong\widehat{A^{**}}\mid_{\widehat{\Delta(A)}}$, then $A$ is a BSE algebra.
 
The BSE- norm property for $A^{**}$ and proving that $A$ and $A^{**}$ are  BSE-norm algebras are assessed in this article.

%$$$$$$$$$$$$$$$$$$$$$$$$$$$$$$$$$$$$$$$$$$$$$$$$$$$$$$$$$$$$$$$$$$$$$$$$$$$$$$$$$$$$$$$
	
	\section{Preliminaries}
	\label{sec:pre}
                           The basic terminologies and the related information on  BSE  algebras are extracted from \cite{E6}, \cite{E7}, and \cite{E8}.

%\subsection{-}
Let $A$ be a  commutative semisimple Banach algebra, and $\Delta(A)$ be the character space of $ A$ with the Gelfand topology. In this study, $\Delta(A)$ represents the set of all non-zero multiplicative
linear functionals over $A$.

 Assume that $C_{b}(\Delta(A))$ is the  space consisting of all complex-valued continuous and bounded functions on $\Delta(A)$.
A linear operator $T$ on $A$ is named a multiplier if for  all $x,y\in A$, $T(xy)=xT(y)$.
The set of all multipliers on $A$ will be expressed as $M(A)$. It is obvious that $M(A)$ is a Banach algebra, and if $A$ is an unital Banach algebra, then 
$M(A)\cong A$.  As observed in   \cite{klar},  for each $T\in M(A)$ there exists a unique bounded  continuous function $\widehat{T}$ on $\Delta(A)$ expressed as:
$$\varphi(Tx)=\widehat{T}(\varphi)\varphi(x),$$
for all $x\in A$ and $\varphi\in \Delta(A)$.
By setting $\{\widehat{T}: T\in M(A)\}$, the $\widehat{M(A)}$ is yield.
%\subsection{-} 

 If the Banach algebra $A$ is  semisimple, then the Gelfand map $\Gamma_{A}: A\to\widehat A$, $f\mapsto \hat f$, is injective or equivalently, and the following equation holds:
\begin{align*}
\underset{\varphi\in\Delta(A)}{\bigcap}\ker(\varphi)=\{0\}.
\end{align*}

A bounded complex-valued continuous function $\sigma$ on $\Delta(A)$  is named a  BSE function, if there exists a positive real number $\beta$ in a sense that for every finite complex-number $c_{1},\cdots,c_{n}$,  and  the same many $\varphi_{1},\cdots,\varphi_{n}$ in $\Delta(A)$ the following inequality
$$\mid\sum_{j=1}^{n}c_{j}\sigma(\varphi_{j})\mid\leq \beta\|\sum_{j=1}^{n}c_{j}\varphi_{j}\|_{A^{*}}$$
holds.\\
The set of all BSE-functions is expressed by $C_{\textup{BSE}}(\Delta( A))$, where for 
 each $\sigma$, the BSE-norm of $\sigma$, $\|\sigma\|_{\textup{BSE}}$ 
is   the infimum of all  $\beta$s  applied in the above inequality.   That $(C_{\textup{BSE}}(\Delta(A)), \|.\|_{\textup{BSE}})$ is a semisimple Banach subalgebra of $C_{b}(\Delta(A))$ is in Lemma 1 proved in \cite{E6}. Algebra $A$ is named a BSE algebra  if it meets the following condition:
$$\widehat{M(A)}\cong C_{\textup{BSE}}(\Delta(A)).$$ 
If $A$ is unital, then $\widehat{M(A)}\cong \widehat{A}\mid_{\Delta(A)}$, indicating that $A$ is a  BSE algebra if and only if  $C_{\textup{BSE}}(\Delta(A))\cong \widehat{A}\mid_{\Delta(A)}$.
The semisimple Banach algebra $A$ is named a norm-BSE algebra if there exists some  $K>0$ in a sense that for each $a\in A$, the following holds: 
$$
\|a\|_{A}\leq K\|\hat{a}\|_{\textup{BSE}}
$$
%**********************************************************************8

Let $A$ be a Banach algebra,  $F,G\in A^{**}$, $f\in A^{*}$ and $x,y\in A$.\\
1) The first Arens multiplication is defined by:
\begin{align*}
F\square G(f) &=F(G.f)\\
G.f(x) &=G(f.x)\\
f.x(y) &= f(xy) 
\end{align*}
2) The second Arens multiplication on $A^{**}$ is defined by:
\begin{align*}
  F\lozenge G(f) &= G(f.F)\\
             f.F(x)   &=F(x.f)\\
             x.f(y)   &=f(yx)
\end{align*}
It is obvious that   $(A^{**},\square)$ and $(A^{**},\lozenge)$  are Banach algebras.
If the binary operation  $\square$ is the same as the action of $\lozenge$, then $A$ is named the  Arens regular Banach algebra.    
If $A$ is a  commutative Banach  algebra, then $A^{**}$ is a commutative Banach  algebra if and only if $A$ is an Arens regular
 Banach algebra; \cite{DL}.
%$$$$$$$$$$$$$$$$$$$$$$$$$$$$$$$$$$$$$$$$$$$$$$$$$$$$$$$$$$$$$$$$$$$$$$$$$$$$$$$$$$$$$$$$$$$$$$$$$$$$$$$$$$$$$$$$$$$$$$$$$$$

	\section{BSE algebras and their second Dual }	
           
       The structure of the  BSE functions on $\Delta(A^{**})$ is characterization and the correlations between the BSE property of $A$ and their second duals are assessed.
 
Let
  $A$ be a commutative semisimple Banach algebra. Then\\
 (i) $C_{\textup{BSE}}(\Delta(A))$ equals the set of all $\sigma\in C_{b}(\Delta(A))$ for which there exists at least one bounded
 net $(a_{\alpha})$ in $A$ with $\lim \widehat{a_{\alpha}}(\varphi)= \sigma(\varphi)$ for all $\varphi\in \Delta(A)$.\\
 (ii) $C_{\textup{BSE}}(\Delta(A))= C_{b}(\Delta(A))\cap A^{**}$, where $A^{**}$ is the second dual of $A$. \\
 Proved in [Theorem4 , \cite{E6}].
  %$$$$$$$$$$$$$$$$$$$$$$$$$$$$$$$$$$$$$$$$$$$$$$$$$$$$$$$$$$$$$$$$$$$$$$$$$
 \begin{lem}\label{lmde}
 Let $A$ be a commutative semisimple  Banach algebra. Then
 $$
 {\overline{<\Delta(A)>}}^{w^{*}}= A^{*}
 $$
 \end{lem}
 \begin{proof}
Assume that $B= {\overline{<\Delta(A)>}}^{w^{*}}$ and $B\subsetneqq A^{*}$. Defne $\Theta: B\longrightarrow\mathbb{C}$  by $\Theta(b)=0$ for all $b\in B$. Let $b_{0}\in A^{*}\setminus B$. Then 
by referring to the Hahn–Banach theorem there exists a $w^{*}$- continuous linear function $\overline{\Theta}: A^{*}\longrightarrow\mathbb{C}$ where $\overline{\Theta}(b_{0})=1$ and $\overline{\Theta}(b)= 0$
for all $b\in B$. Because  $(A^{*},w^{*})$ is a locally convex and $\overline{\Theta}$ is $w^{*}$-continuous,
thus there exists $a\in A$ where $\overline{\Theta}= \hat{a}$. If $\varphi\in \Delta(A)$, the following is yield:
\begin{align*}
\varphi(a)= \hat{a}(\varphi)=\overline{\Theta}(\varphi)=0 
\end{align*}
Because $A$ is semisimple,   $a=0$. Thus,
$$
\overline{\Theta}(b_{0})= \hat{a}(b_{0})= b_{0}(0)= 0
$$
which is a contradiction; therefore $ {\overline{<\Delta(A)>}}^{w^{*}}= A^{*}$.
 \end{proof}
%$$$$$$$$$$$$$$$$$$$$$$$$$$$$$$$$$$$$$$$$$$$$$$$$$$$$$$$$$$$$$$$$$$$$$
 \begin{lem}\label{lD}
 Let $A$ be a commutative Banach algebra and $ {\overline{<\Delta(A)>}}^{w^{*}}= A^{*}$. Then 
 $\Delta(A^{**})= {\overline{\widehat{\Delta(A)}}}^{w^{*}}$.
 \end{lem}
\begin{proof}
Provided in \cite{DL}.
\end{proof}
%**********************************************************************
 \begin{cor}\label{lsim}
 Let $A$ be a commutative  semisimple  unital Banach algebra. Then 
 $$
 \Delta(A^{**})=\widehat{\Delta(A)}
 $$
 \end{cor}
 \begin{proof}
  Because $A$ is unital, $\Delta(A)$ is 
$w^{*}$- compact, thus, $\widehat{\Delta(A)}= {\overline{\widehat{\Delta(A)}}}^{w^{*}}$.
Then by applying Lemma \ref{lmde} and Lemma \ref{lD}, the equation $\Delta(A^{**})=\widehat{\Delta(A)}$ is yield.
 \end{proof}
 %@@@@@@@@@@@@@@@@@@@@@@@@@@@@@@@@@@@@@@@@@@@@@
 \begin{rem}
  Let $A$ be a commutative  semisimple  unital Banach algebra. Assume that $\varphi_{1},\cdots,\varphi_{n}\in\Delta(A)$ and $c_{1},\cdots, c_{n}\in\mathbb{C}$, then the following is yield:
  $$
  \|\sum_{i=1}^{n}c_{i}\varphi_{i}\|_{A^{*}}= \|\sum_{i=1}^{n}c_{i}\widehat{\varphi_{i}}\|_{A^{***}}
  $$
 \end{rem}

\begin{rem}
Let $A$ be a commutative    Arens regular   Banach algebra. If $A^{**}$ is semisimple, then
$A$ is semisimple. Because $\widehat{\ker(\varphi)} \subseteq \ker(\hat{\varphi})$ and  $\widehat{\Rad(A)}\subseteq \Rad(A^{**})$. then $A$ is semisimple, because $A^{**}$ is semisimple, .
\end{rem}

\begin{thm}
Let $A$ be a commutative  Arens regular semisimple unital  Banach algebra. Then the following statement holds:
$$
C_{\textup{BSE}}(\Delta(A^{**}))\cong C_{\textup{BSE}}(\Delta(A))
$$
These two as Banach algebras are isometric;
\end{thm}
\begin{proof}
Because $A$ is semisimple, by applying Corollary \ref{lsim}, the $\Delta(A^{**})=\widehat{\Delta(A)}$ is yield, thus, 
$C_{b}(\Delta(A^{**}))\cong C_{b}(\Delta(A))$, and $A^{****}\mid_{\widehat{\Delta(A)}}\cong A^{**}\mid_{\Delta(A)}$.
 Consequently by applying \cite[Theorem 4]{E6}, $C_{\textup{BSE}}(\Delta(A^{**}))\cong C_{\textup{BSE}}(\Delta(A))$ is yield.
\end{proof}

The correlation between the BSE-property of  algebra $A$ and it's second dual $A^{**}$ is assessed as follows:
\begin{lem}\label{la1}
Let  $A$ be a commutative unital Banach algebra such that $A$ is Arens regular and $A^{**}$ is semisimple. If $A$ is a BSE algebra, so is $A^{**}$.
\end{lem}
%$$$$$$$$$$$$$$$$$$$$$$$$$$$$$$$$$$$
\begin{proof}
Let $\Sigma\in C_{\textup{BSE}}(\Delta(A^{**}))$. Define $\sigma(\varphi)= \Sigma(\widehat{\varphi})$ for all $\varphi\in\Delta(A)$. Obviously, $\sigma$ is well defined and $\sigma\in C_{\textup{BSE}}(\Delta(A))$. By allowing $\varphi_{1},\cdots, \varphi_{n}\in \Delta(A)$ and $c_{1}, \cdots, c_{n}\in\mathbb{C}$
the following is yield:
\begin{align*}
\mid\sum_{j=1}^{n}c_{j}\sigma(\varphi_{j})\mid &= \mid\sum_{j=1}^{n}c_{j}\Sigma(\widehat{\varphi_{j}})\mid\\
                                                                          &\leq M \mid\mid\sum_{j=1}^{n}c_{j}\widehat{\varphi_{j}}\mid\mid_{A^{***}}\\
                                                                          &= M\mid\mid\sum_{j=1}^{n}c_{j}{\varphi_{j}}\mid\mid_{A^{*}}
\end{align*}
where $M>0$. Because $C_{\textup{BSE}}(\Delta(A))=\widehat{A}\mid_{\Delta(A)}$, there exists $a\in A$ where $\sigma = \hat{a}$. \\
Assume that $\Phi\in \Delta(A^{**})$, thus $\Phi=\hat{\varphi}$ for some  $\varphi\in\Delta(A)$. Which yield:
\begin{align*}
\Sigma(\Phi)= \Sigma(\hat{\varphi}) &=\sigma(\varphi)\\
                                                            &= \hat{a}(\varphi)\\
                                                              &=\widehat{\hat{a}}(\hat{\varphi})\\
                                                                &= \widehat{\hat{a}}(\Phi)
\end{align*}
 where, $\Sigma= \widehat{\hat{a}}$, therefore $C_{\textup{BSE}}(\Delta(A^{**}))=\widehat{A^{**}}\mid_{\widehat{\Delta(A)}}$.
\end{proof}

\begin{lem}
Let  $A$ be a commutative  unital Banach algebra where $A$ is Arens regular and $A^{**}$ is semisimple. If $A^{**}$ is a BSE algebra and $\widehat{A}\mid_{\Delta(A)}\cong\widehat{A^{**}}\mid_{\widehat{\Delta(A)}}$, then $A$ is a BSE-algebra.
\end{lem}
\begin{proof}
Assume that $A^{**}$ is a BSE algebra, the following is yield:
\begin{align*}
C_{BSE}(\Delta(A)) & \cong C_{BSE}(\Delta(A^{**})\\
                               &  \cong \widehat{A^{**}}\mid_{\widehat{\Delta(A)}}\\
                               &  \cong \widehat{A}\mid_{\Delta(A)}
\end{align*}
thus, $A$ is a BSE- algebra.
\end{proof}
%+++++++++++++++++++++++++++++++++++++++++++++
At this stage, based on the established prerequisite the primary Theorem, is expressed as follows:
\begin{thm}
Let  $A$ be a commutative  unial Banach algebra, such that $A$ is Arens regular and $A^{**}$ is semisimple. 
Then the following statements are equivalent:\\
(i) $A$ is a BSE algebra.\\
(ii) $A^{**}$ is a BSE algebra and 
$\widehat{A}\mid_{\Delta(A)}\cong\widehat{A^{**}}\mid_{\widehat{\Delta(A)}}$.
\end{thm}
%$$$$$$$$$$$$$$$$$$$$$$$$$$$$$$$$$$$$$$$$$$$$$$$$$$
\begin{ex}
Let $(K, d)$ be a compact metric space and $A$ be a unital commutative semisimple Banach algebra. Then $A$ is a BSE-algebra if and only if $Lip(K, A)$ is a BSE-algebra;   \cite{AK}. In this context, the following statements are equivalent:\\
(i) $A$ is a BSE algebra.\\
(ii) $A^{**}$ is a BSE algebra and 
$\widehat{A}\mid_{\Delta(A)}\cong\widehat{A^{**}}\mid_{\widehat{\Delta(A)}}$.\\
(iii) $\Lip(K, A)$ is a BSE-algebra.\\
(iv) $\Lip(K, A^{**})$ is a BSE-algebra and 
$\widehat{A}\mid_{\Delta(A)}\cong\widehat{A^{**}}\mid_{\widehat{\Delta(A)}}$.\\
\end{ex}

%@@@@@@@@@@@@@@@@@@@@@@@@@@@@@@@@@@@@@@@
\begin{ex}
 Let $X$ be a metric space and $A = \Lip (X)$. Then $A$ is unital and BSE algebra;   \cite{AK2}. By applying Lemma \ref{la1}, ${\Lip (X)}^{**}$  becomes a BSE algebra.

\end{ex}
%++++++++++++++++++++++++++++++++++++++++++
\begin{ex}\label{ex1}
Let $A$ be a commutative reflexive Banach algebra. Then 
$$
C_{\textup{BSE}}(\Delta(A))= \widehat{A}\mid_{\Delta(A)}
$$
Moreover, $A$ is unital, then $A$ is a BSE algebra.
\end{ex}
\begin{proof}
According to [Theorem 4, \cite{E6}] because $A$ is reflexive, and the following is yield: 
\begin{align*}
C_{\textup{BSE}}(\Delta(A)) &=A^{**}\mid_{\Delta(A)}\cap C_{b}(\Delta(A))\\
                              &=\widehat{A}\mid_{\Delta(A)}\cap C_{b}(\Delta(A))\\
                              &=\widehat{A}\mid_{\Delta(A)}
\end{align*}
\end{proof}

\section{BSE-norm algebra}
The BSE-norm property of the second dual of Banach algebra is assessed and the correlation between the BSE-norm of the  Banach algebra $A$  and its second dual is assessed.
%$$$$$$$$$$$$$$$$$$$$$$$$$$$$$$$$$$$$$$$$$$$$$$$$$$$$$$$$$$$$$$$$$$$$$$$$$
\begin{thm}
Let  $A$ be a commutative  unial Banach algebra such that $A$ is Arens regular and $A^{**}$ is semisimple. Then $A^{**}$ is a BSE-norm algebra.
\end{thm}
\begin{proof}
By applying the open mapping theorem
there exists some $M>0$ where $||F||_{A^{**}}\leq M||F\mid_{\Delta(A)}||_{A^{**}}$ for each $F\in A^{**}$,  because the map $\Theta: ~ A^{**}\longrightarrow C_{b}(\Delta(A))$ given by $\Theta(F)= F\mid_{\Delta(A)}$ is continuous, linear and injective, thus:
\begin{align*}
\|\hat{F}\|_{\textup{BSE}} &= \sup \{\mid\sum_{i=1}^{n}c_{i}\hat{F}(\Phi_{i})\mid~\Phi_{i}\in\Delta(A^{**}), ~ \|\sum_{i=1}^{n}c_{i}\Phi_{i}\|_{A^{***}}\leq 1\}\\
                           &= \sup \{\mid\sum_{i=1}^{n}c_{i}\hat{F}(\hat{\varphi_{i}})\mid~ \varphi_{i}\in\Delta(A), \|\sum_{i=1}^{n}c_{i}\varphi_{i}\|_{A^{*}}\leq 1\}\\
                            &= \sup \{\mid\sum_{i=1}^{n}c_{i}{F}({\varphi_{i}})\mid~ \varphi_{i}\in\Delta(A), \|\sum_{i=1}^{n}c_{i}\varphi_{i}\|_{A^{*}}\leq 1\}\\
                            &= \|F\mid_{\Delta(A)}\|_{\textup{BSE}}
\end{align*}
in addition:
\begin{align*}
\|F\mid_{\Delta(A)}\|_{A^{**}} &= \sup\{\mid F(\varphi)\mid :\varphi\in\Delta(A)\}\\
                                                &\leq \|F\mid_{\Delta(A)}\|_{\textup{BSE}}
\end{align*}
Then by applying the open mapping theorem, the following is yield:
$$
\|F\|_{A^{**}}\leq M \|F\mid_{\Delta(A)}\|_{A^{**}}\leq M\|F\mid_{\Delta(A)}\|_{\textup{BSE}}=  M \|\hat{F}\|_{\textup{BSE}}
$$
therefore, $A^{**}$ is a BSE-norm algebra.
\end{proof}

%_______________________________________________________
\begin{thm}
Let  $A$ be a commutative unital Banach algebra, such that $A$ is Arens regular
and $A^{**}$ is semisimple. Then $A$ is a BSE-norm algebra.
\end{thm}
\begin{proof}
Because $A^{**}$ is a BSE-norm algebra, consequently, by definition there exists $M>0$ where $\|F\|_{A^{**}}\leq M{\|\widehat{F}\|_{\textup{BSE}}}$ for each $F\in A^{**}$. If $a\in A$, the following is yield: 
$$
\|a\|_{A} = \|\hat{a}\|_{A^{**}}\leq M\|\widehat{\hat{a}}\|_{\textup{BSE}} 
$$
in addition
\begin{align*}
\|\widehat{\hat{a}}\|_{\textup{BSE}} &= \sup\{\mid\sum_{i=1}^{n}c_{i}\Phi_{i}(\hat{a})\mid ~\mid~ \Phi_{i}\in\Delta(A^{**}), \|\sum_{i=1}^{n}c_{i}\Phi_{i}\|_{A^{***}}\leq 1\}\\
                                           &= \sup\{\mid\sum_{i=1}^{n}c_{i}\widehat{\varphi_{i}}(\hat{a})\mid ~\mid~ \|\sum_{i=1}^{n}c_{i}\widehat{\varphi_{i}}\|_{A^{***}}\leq 1\}\\
                                         &=  \sup\{\mid\sum_{i=1}^{n}c_{i}\varphi_{i}(a)\mid ~\mid~ \|\sum_{i=1}^{n}c_{i}\varphi_{i}||_{A^{*}}\leq 1\}\\ 
                                         &= \|\hat{a}\|_{\textup{BSE}}
\end{align*}
therefore, $\|a\|_{A}\leq M\|{\hat{a}}\|_{\textup{BSE}} $. Then $A$ is a  BSE-norm algebra.
\end{proof}
%^^^^^^^^^^^^^^^^^^^^^^^^^^^^^^^^^^^^^
\begin{cor}
Let  $A$ be a commutative unital Banach algebra, such that $A$ is Arens regular
and $A^{**}$ is semisimple. Then $A$ and $A^{**}$ are  BSE-norm algebras.
\end{cor}
%$$$$$$$$$$$$$$$$$$$$$$$$$$$$$$$$$$$$$$$$$$$$$$$$$$$$$$$$$$$$
\begin{ex}\label{lc}
Let $A$ be a commutative semisimple reflexive Banach algebra where $\Delta(A)$ is a discrete space.
Then $A$ is a BSE-norm algebra.
\end{ex}
\begin{proof} 
By applying Example \ref{ex1},  $C_{\textup{BSE}}(\Delta(A))=  {\hat{A}}\mid_{\Delta(A)}$,
making ${\hat{A}}\mid_{\Delta(A)}$ is a closed BSE-norm space, hence $\Phi:~ A\longrightarrow {\hat{A}}\mid_{\Delta(A)}$ is an one-to-one  continuous map where $\Phi(A)= {\hat{A}}\mid_{\Delta(A)}$. Then by applying the 
open mapping theorem it is revealed that there exists some $M>0$ where $\|a\|_{A}\leq M \|\hat{a}\|_{\textup{BSE}}$ for each $a\in A$,
therefore $A$ is a BSE-norm algebra.
\end{proof}
\begin{ex}
1) If $A$ is a commutative unital finite  dimensional Banach algebra , then $A$ is a BSE-norm algebra.\\
2) Let $G$ be a compact group and $A= l_{p}(G)$. Then   $A$ is a BSE-norm algebra.\\
\end{ex}
\begin{proof}
1)
Since $A$ is finite  dimensional,  $A\cong C^{n}$, for some positive integer $n$ under equivalent norms.    Example \ref{lc} implies that $A$ is a BSE-norm algebra.\\
2) Because $\Delta(A)= \widehat{G}$ is a discrete space,  $A$  is a semisimple and reflexive  algebra, thus by according to Example \ref{lc}, $A$ is a BSE-norm algebra.
\end{proof}
\begin{ex}
Assume that $G$ is a compact group and $A$ is a reflexive  semisimple Banach algebra where $\Delta(A)$ is discrete
and $1<p<\infty$. Then $L_{p}(G,A)$ is a BSE-norm algebra.
\end{ex}
\begin{proof}
The following can be written:
\begin{align*}
L_{p}(G,A)^{**} &\cong L_{q}(G,A^{*})^{*}\\
                           &\cong L_{p}(G,A^{**})\cong L_{p}(G,A)
\end{align*}
Consequently  $L_{p}(G,A)$ is reflexive, and  
$$
\Delta(L_{p}(G,A))= \widehat{G}\times \Delta(A)
$$
Because $G$ is a compact group, and $\Delta(A)$ is  discrete, it is concluded  that $\Delta(L_{p}(G,A))$ is  discrete.
Therefore by applying Example \ref{lc}, $L_{p}(G,A)$ is a BSE-norm algebra.
\end{proof}

\section{Acknowledgment}
           Appreciation is extended to Assistant  Professor  Hojjat Farzadfard for his valuable comments and suggestions regarding the manuscript.

\vspace{9mm}

{\footnotesize \noindent
	M. Aghakoochaki\\
	Department of Pure Mathematis\\
	Faculty of Mathematics and Statistics\\
	University of Isfahan\\
	Isfahan, 81746-73441\\
	Iran\\
mkoochaki@sci.ui.ac.ir\\
	
	\noindent
	A. Rejali\\
	Department of Pure Mathematis\\
	Faculty of Mathematics and Statistics\\
	University of Isfahan\\
	Isfahan, 81746-73441\\
	Iran\\
	rejali@sci.ui.ac.ir\\

\end{document}